\documentclass[11pt]{article}

\usepackage{epsfig}
\usepackage{amssymb}

\newtheorem{theorem}{Theorem}[section]

\newtheorem{proposition}[theorem]{Proposition}
\newtheorem{definition}[theorem]{Definition}

\newtheorem{algorithm}[theorem]{Algorithm}

\newcommand{\cc}{{\mathbb C}}
\newcommand{\rr}{{\mathbb R}}
\newcommand{\zz}{{\mathbb Z}}
\newcommand{\ff}{{\mathbb F}}
\newcommand{\nn}{{\mathbb N}}

\newcommand{\bfa}{{\mathbf a}}

\newcommand{\bfp}{{\mathbf p}}
\newcommand{\bfq}{{\mathbf q}}

\newcommand{\bfv}{{\mathbf v}}
\newcommand{\bfx}{{\mathbf x}}
\newcommand{\bfz}{{\mathbf z}}

\newcommand\qed{{\hspace*{\fill}$\Box$\vskip12pt plus 1pt}}

\begin{document}

\title{Tropical Algebraic Geometry in Maple \\
a preprocessing algorithm for finding common factors
to multivariate polynomials with approximate 
coefficients\thanks{{\em Date:} 1 September 2008.
This material is based upon work
supported by the National Science Foundation
under Grant No.\ 0713018.}
}

\author{
Danko Adrovic\thanks{
Department of Mathematics, Statistics, and Computer Science,
University of Illinois at Chicago, 851 South Morgan (M/C 249),
Chicago, IL 60607-7045, USA.
{\em Email:} adrovic@math.uic.edu}
\and
Jan Verschelde\thanks{
Department of Mathematics, Statistics, and Computer Science,
University of Illinois at Chicago, 851 South Morgan (M/C 249),
Chicago, IL 60607-7045, USA.
{\em Email:} jan@math.uic.edu or jan.verschelde@na-net.ornl.gov
{\em URL:} http://www.math.uic.edu/{\~{}}jan}
}

\date{{\em to Keith Geddes, on his $60^{th}$ birthday}}

\maketitle

\begin{abstract}
Finding a common factor of two multivariate polynomials with
approximate coefficients is a problem in symbolic-numeric computing.
Taking a tropical view on this problem leads to efficient
preprocessing techniques, applying polyhedral
methods on the exact exponents with numerical techniques on
the approximate coefficients.  With Maple we will illustrate
our use of tropical algebraic geometry.
\end{abstract}

\section{Introduction}

Tropical algebraic geometry is a relatively new language
to study skeletons of algebraic varieties.
Introductions to tropical algebraic geometry
are in~\cite{RST05} and~\cite[Chapter~9]{Stu02}.
Computational aspects are addressed in~\cite{BJSST07} and~\cite{The06}.  
One goal of this paper is to explain some new words of this language, 
and to show how a general purpose computer algebra system like Maple
is useful to explore and illustrate tropical algebraic geometry.
For software dedicated to tropical geometry, we refer
to Gfan~\cite{Jen07,Jen08}, a SINGULAR library~\cite{JMM07},
and TrIm~\cite{SY08}.

The roots of tropical algebraic geometry run as deep
as the work of Puiseux~\cite{Pui1850} and Ostrowski~\cite{Ost22},
therefore our focus is
on answering a practical question in computer algebra:
{\em when do two polynomials have a common factor?}
Viewing this question in tropical algebraic geometry leads
to a symbolic-numeric algorithm.  In particular, we will say
that {\em tropisms give the germs to grow the tentacles of the
common amoeba.}  The paper is structured in four parts, each
part explaining one of the key concepts of the tropical sentence.

Our perspective on tropical algebraic geometry originates
from polyhedral homotopies~\cite{HS95}, \cite{Li03}, \cite{VVC94}
to solve polynomial systems implementing Bernshte\v{\i}n's
first theorem~\cite{Ber75}.
Another related approach that led to tropical mathematics
is idempotent analysis.  In~\cite{Reg08}, a Maple package
is presented for a tropical calculus with application to
differential boundary value problems.

Related work on our problem concerns the factorization
of sparse polynomials via Newton polytopes~\cite{Abu08}, \cite{GL01};
approximate factorization~\cite{CGHKW01}, \cite{CGKW02}, \cite{GvH07},
\cite{KMYZ08}, \cite{SVW8},
and the GCD of polynomials with approximate coefficients~\cite{ZD04}.
The polynomial absolute factorization is also addressed in~\cite{Che05}
and the lectures in~\cite{CG05} offer a very good overview.
Criteria based on polytopes for the irreducibility of polynomials
date back Ostrowski~\cite{Ost22}.
In this paper we restrict our examples to polynomials in two
variables and refer to polygons instead of polytopes.
The terminology extends to general dimensions and polytopes,
see~\cite{Zie95}.

That two polynomials with approximate coefficients have a common
factor is quite an exceptional situation.  Therefore it is important
to have efficient preprocessing criteria to decide quickly.
The preprocessing method we develop in this paper attempts to build
a Puiseux expansion starting at a common root at infinity.
To determine whether a root at infinity is isolated or not
we apply the Newton-Puiseux method, extending the proof
outlined by Robert Walker in~\cite{Wal50}, see also~\cite{dJP00},
towards Joseph Maurer's general method~\cite{Mau80} for space curves.  
A more algorithmic method than~\cite{Mau80} is given in~\cite{AMNR92}
along with an implementation in CoCoA.
General fractional power series solutions are described in~\cite{McD02}.
See~\cite{JMM08}, \cite{JMSW08} and~\cite{PR08}
for recent symbolic algorithms, and~\cite{Pot07}, \cite{PR08b} for
a symbolic-numeric approach.  
The complexity for computing Puiseux expansions for plane curves
is polynomial~\cite{Wal00} in the degrees.
As an alternative to Puiseux series, extended Hensel series have
been developed, see~\cite{Sas08} for a survey.  Good numerical
convergence is reported in~\cite{IS07}.

We show that via suitable coordinate transformations,
the problem of deciding whether there is a common factor
is reduced to univariate root finding, with the univariate
polynomials supported on edges of the Newton polygons of
the given equations.  Also in the computation of the second
term of the Puiseux series expansion, we do not need to utilize
all coefficients of the given polynomials.  In the worst case,
the cost of deciding whether is a factor is a cubic polynomial 
in the number of monomials of the given polynomials.

Certificates for the existence of a common factor consist
of exact and approximate data: the exponents and coefficients
of the first two terms of a Puiseux series expansion of the
factor at a common root at infinity.  The leading exponents
of the Puiseux expansion form a so-called {\em tropism}~\cite{Mau80}.
The coefficients are numerical solutions of overdetermined systems.
If a more explicit form of the common factor is required,
more terms in the Puiseux expansion can be computed up to precision 
needed for the application of sparse interpolation techniques,
see~\cite{CL08}, \cite{GLL06}, \cite{KL03}, and~\cite{Lee07}.

The {\tt ConvexHull} and {\tt subs} commands of Maple are very valuable
in implementing an interactive prototype of the preprocessing algorithm.
For explaining the intuition behind the algorithm, we start illustrating
amoebas and the tentacles.  Once we provide an abstraction for the tentacles 
we give an outline of the algorithm and sketch its cost.

\noindent {\bf Acknowledgements.}
This paper is based on the talk the second author gave
at MICA 2008 -- Milestones in Computer Algebra, a conference in
honour of Keith Geddes' 60th Birthday,
Stonehaven Bay, Trinidad and Tobago, 1-3 May.
We thank the organizers for this wonderful and inspiring conference
for their invitation.

\section{Amoebas}

Looking at the asymptotics of varieties gives a natural
explanation for the Newton polygon.  This polygon will provide
a first classification of the approximate coefficients of the
given polynomials.  This means that at first we may ignore
coefficients of monomials whose exponents lie in the interior
of the Newton polygon.

\subsection{asymptotics of varieties}

Our input data are polynomials in two variables $x$ and $y$.
The set of values for $x$ and $y$ that make the polynomials zero
is called a variety.  
Varieties are the main objects in algebraic geometry.  
In 1971, G.M. Bergman~\cite{Ber71} considered logarithms of varieties.
In tropical algebraic geometry, we look at the asymptotics
of varieties.  
\begin{equation}
  \begin{array}{ccccc}
     \log & : & \cc^* \times \cc^* & \rightarrow & \rr \times \rr \\
          &   & (x,y) & \mapsto & (\log(|x|),\log(|y|))
  \end{array}
\end{equation}
Because the logarithm is undefined at zero, we exclude the origin
restricting the domain of our polynomials to the torus $(\cc^*)^2$,
$\cc^* = \cc \setminus \{ 0 \}$.  Following~\cite{GKZ94}, we arrive
at our first new word~\cite{Vir02}.

\begin{definition}[Gel'fand, Kapranov, and Zelevinsky 1994] {\rm
The {\em amoeba} of a variety is its image under the log map. }
\end{definition}

\begin{figure}[hbt]
\begin{center}
\begin{picture}(400,150)(0,0)
\put(0,0){\epsfig{figure=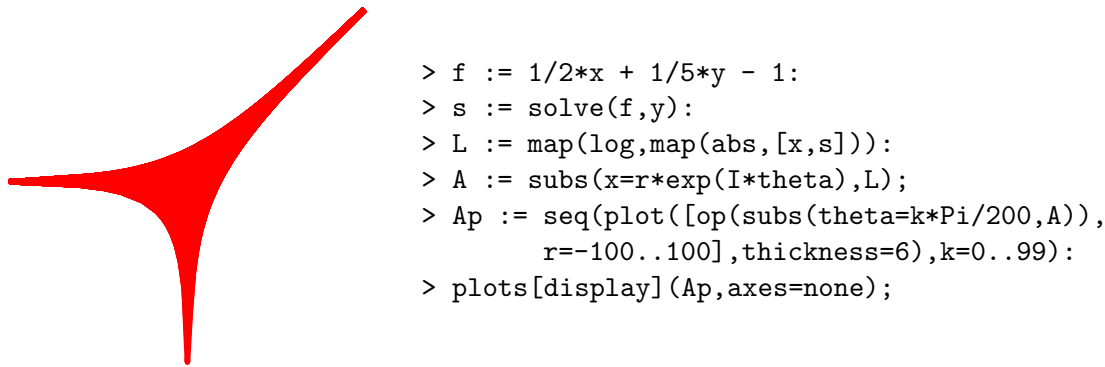,width = 5cm}}
\put(150,70){
\begin{tabular}{l}
{\tt > f := 1/2*x + 1/5*y - 1:} \\
{\tt > s := solve(f,y):} \\
{\tt > L := map(log,map(abs,[x,s])):} \\
{\tt > A := subs(x=r*exp(I*theta),L);} \\
{\tt > Ap := seq(plot([op(subs(theta=k*Pi/200,A)),} \\
{\tt ~~~~~~~ r=-100..100],thickness=6),k=0..99): } \\
{\tt > plots[display](Ap,axes=none);}
\end{tabular}
}
\end{picture}
\caption{The amoeba of a linear polynomial, with all Maple commands
at the right.}
\label{figamoeba1}
\end{center}
\end{figure}

To see what amoebas look like, we use the plotting capabilities
of Maple.  We use polar coordinates (see~\cite{The02} for more
on plotting amoebas) to express a linear variety:

\begin{equation}
f := \frac{1}{2} x + \frac{1}{5} y - 1 = 0
\quad \quad
A := \left[ 
  \ln\left( \left| r e^{I \theta} \right| \right) ,
  \ln\left( \left| \frac{5}{2} r e^{I \theta} - 5 \right| \right)
\right].
\end{equation}
In Figure~\ref{figamoeba1} we see the result of a Maple plot.

\subsection{compactifying amoebas lead to Newton polytopes}

We compactify the amoeba of $f^{-1}(0)$ by taking lines perpendicular
to the tentacles.  As each line cuts the plane in half, we keep those
halves of the plane where the amoeba lives.  The intersection of all
half planes defines a polygon.

\begin{figure}[hbt]
\begin{center}
\begin{picture}(250,140)(0,0)
\put(0,0){\epsfig{figure=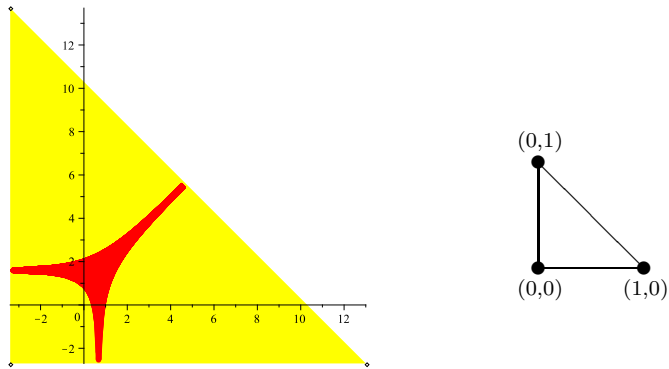,width = 5cm}}
\put(200,40){
\begin{picture}(50,50)(0,0)
\put(0,0){\circle*{5}}   \put(-8,-8){${}_{(0,0)}$}
\put(40,0){\circle*{5}}  \put(32,-8){${}_{(1,0)}$}
\put(0,40){\circle*{5}}  \put(-8,48){${}_{(0,1)}$}
\put(0,0){\line(1,0){40}}
\put(0,0){\line(0,1){40}}
\put(0,40){\line(1,-1){40}}
\end{picture}
}
\end{picture}
\caption{The compactification of the amoeba: the edges of the Newton
polygon (displayed at the right) are perpendicular to the tentacles 
of the amoeba.}
\label{figamoeba2}
\end{center}
\end{figure}

The resulting polygon is the Newton polygon of~$f$.
For the amoeba in Figure~\ref{figamoeba1},
its compactification is shown in Figure~\ref{figamoeba2}.
On Figure~\ref{figamoeba2} we recognize the shape of the triangle,
the Newton polygon of a linear polynomial.  In general, there is a
map~\cite{Sot03} that sends every point in the variety to the interior
of the Newton polygon of the defining polynomial equation.

This geometric derivation of the Newton polygon coincides
with the more formal definition.

\begin{definition} {\rm  
For ${\displaystyle f(x,y) = \sum_{(i,j) \in A} c_{i,j} x^i y^j}$, 
$c_{i,j} \in \cc^*$.
$A$ is the {\em support} of $f$.
The convex hull of $A$ is the {\em Newton polygon}.}
\end{definition}

The Newton polygon models the sparse structure of a polynomial.
Most polynomials arising in practical applications have few monomials
with nonzero coefficients and are called {\em sparse}.  
The Newton polygon assigns
additional significance to the coefficients.  Coefficients associated
to monomials whose exponents span a vertex of the Newton polygon are
more important than coefficients whose exponents lie in the interior
of the Newton polygon.

Plotting amoebas is actually computationally quite involved
-- the use of homotopy continuation methods~\cite{Ver99}
is suggested in~\cite{The02}.
See~\cite{Nis08} for more about amoebas.
We will see that the asymptotics of the amoebas will lead
to a natural reduction of our problem to smaller polynomials
{\em in one variable}.

\section{Tentacles}

The tentacles of the amoeba stretch out to infinity
and are represented by the inner normals, perpendicular
to the edges of the Newton polygon.

\subsection{the direction of the tentacles towards infinity}

Our problem may be stated as follows:
Given two polynomials in two variables 
with {\em approximate} complex coefficients, is there a common factor?

Looking at the problem from a tropical point of view,
we first have the amoeba of the common factor in mind.
The tentacles of that amoeba stretch out to infinity
and these tentacles are represented by inner normals,
perpendicular to the lines at infinity corresponding
to the edges of the Newton polygon.
Because the factor is common to both polynomials,
the normals must be common to both polygons. 
So if there is a factor, there must be at least one
pair of edges with the same inner normal vector.
Such inner normal vector is a {\em tropism}, defined below.

The tropical view will lead to solving the problem first
at infinity, providing an efficient preprocessing criterion.
We first formalize the geometric idea in a proposition.
\begin{proposition} \label{propamoebas}
Let $f$ and $g$ be two polynomials.
If the amoebas of $f$ and $g$ have no tentacle stretching
out to infinity in the same direction, 
then $f$ and $g$ have no common factor.
\end{proposition}
Verifying the conditions of Proposition~\ref{propamoebas}
seems nontrivial at first.

\subsection{normal fans and tropicalization}

\begin{figure}[hbt]
\begin{center}
\begin{picture}(300,150)(0,0)
\put(0,0){\epsfig{figure=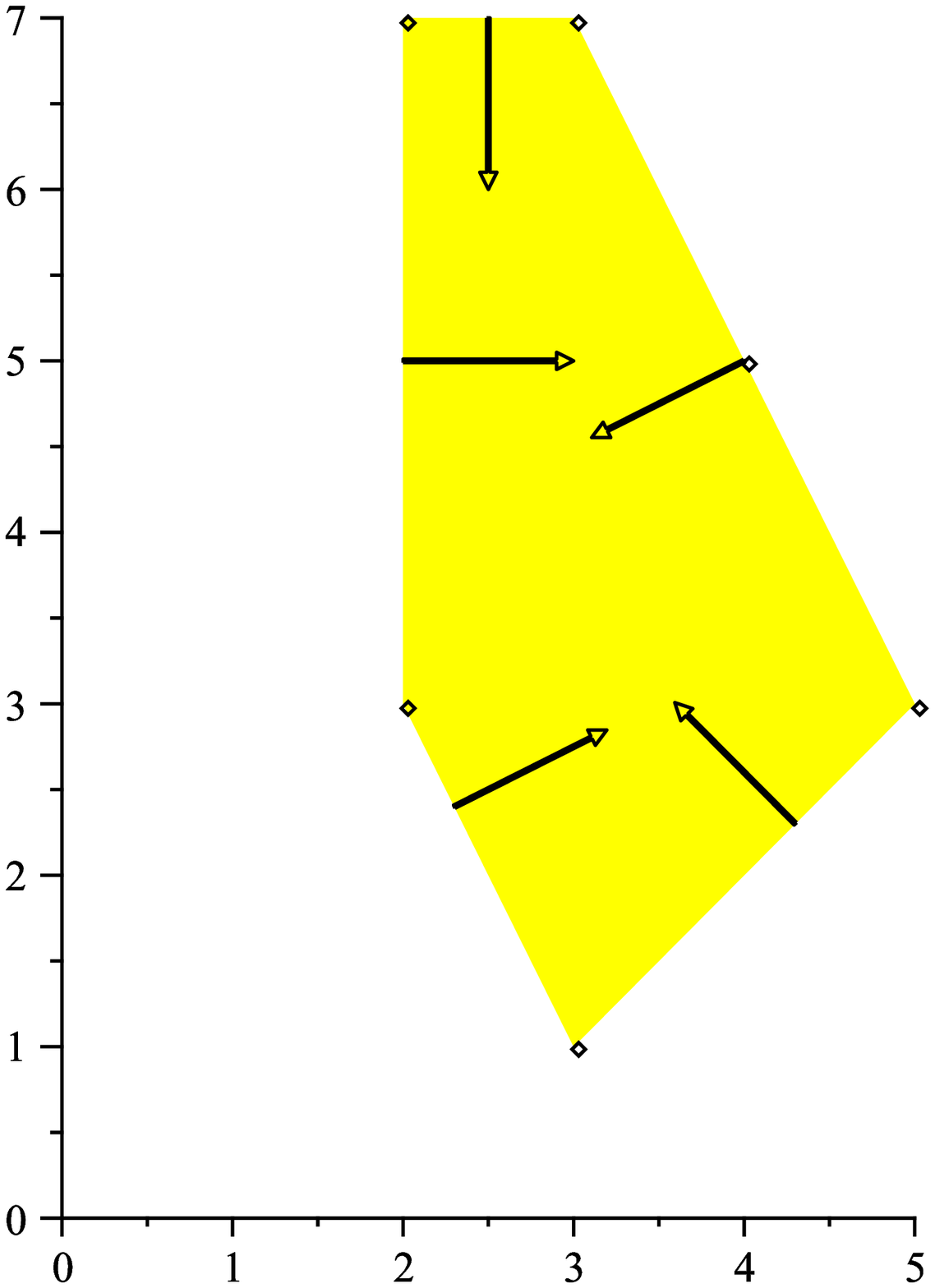,width = 5cm}}
\put(150,0){\epsfig{figure=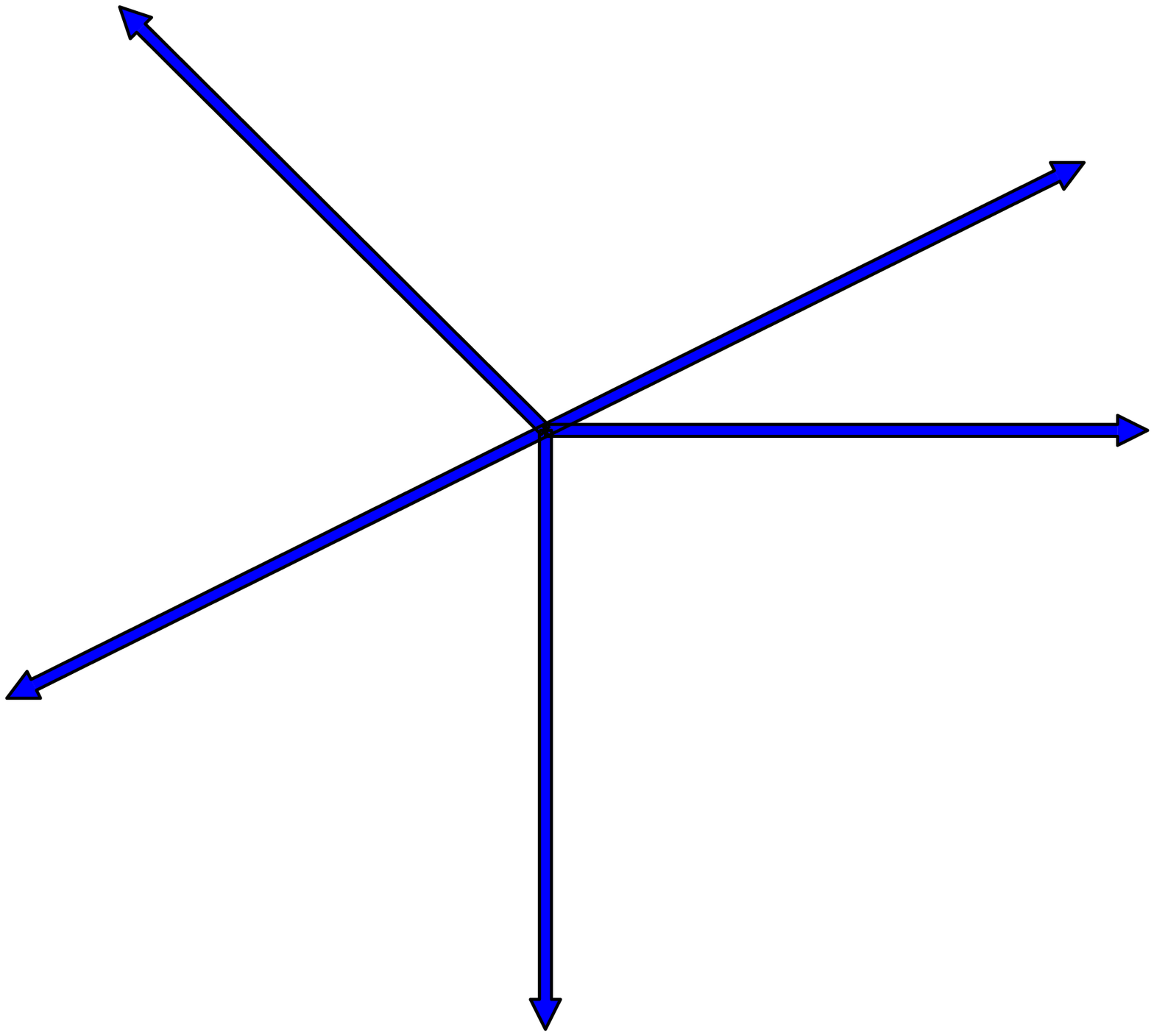,width = 5cm}}
\end{picture}
\end{center}
\caption{The Newton polygon and its normal fan.}
\label{fignormfan}
\end{figure}

The inward pointing normal vectors to the edges represent
the tentacles of the amoeba.
Consider for example

\begin{equation}
f := x^3 y + x^2 y^3 + x^5 y^3 + x^4 y^5 + x^2 y^7 + x^3 y^7.
\end{equation}
In Figure~\ref{fignormfan} we show the Newton polygon of~$f$
and its normal fan.

The collection of inner normals to the edges of the Newton polygon
forms {\em a tropicalization} of~$f$, denoted by~${\rm Trop}(f)$.
To formalize this notion, we introduce the following definitions.

Exponents and direction vectors are related through duality
via the inner product.

\begin{definition} {\rm The {\em inner product} is
\begin{equation}
 \begin{array}{ccccc}
  \langle \cdot , \cdot \rangle & \zz^2 \times \zz^2 & \rightarrow & \zz \\
     & ((i,j),(u,v)) & \mapsto & i u + j v.
 \end{array}
\end{equation} }
\end{definition}

\begin{figure}[hbt]
\begin{center}
\begin{picture}(300,145)(0,0)
\put(-20,0){\epsfig{figure=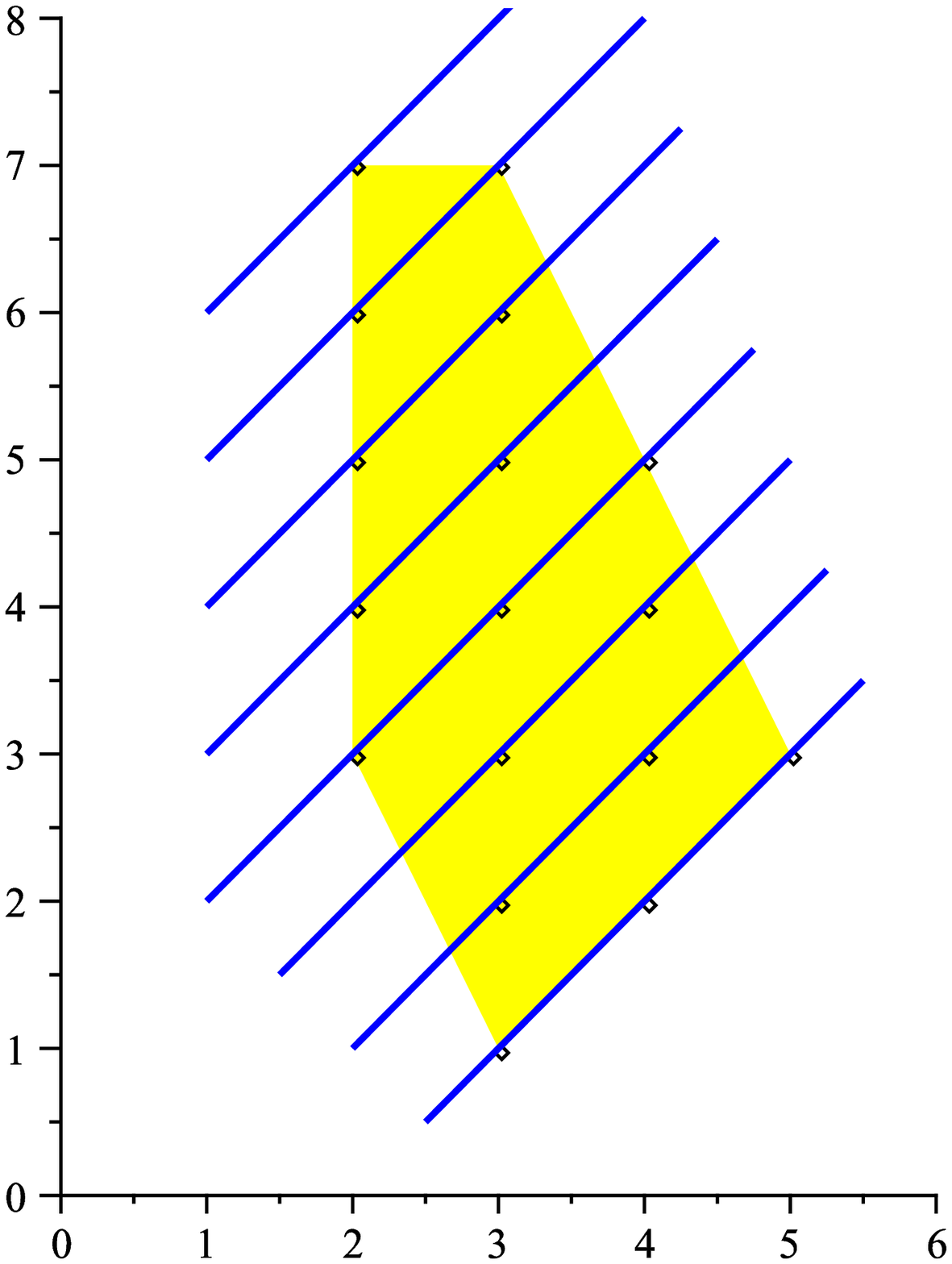,width = 6cm, height = 5cm}}
\put(150,115){$-1 \times i + (+1) \times j = +5$}
\put(150,100){$-1 \times i + (+1) \times j = +4$}
\put(150,85){$-1 \times i + (+1) \times j = +3$}
\put(150,70){$-1 \times i + (+1) \times j = +2$}
\put(150,55){$-1 \times i + (+1) \times j = +1$}
\put(150,40){$-1 \times i + (+1) \times j = 0$}
\put(150,25){$-1 \times i + (+1) \times j = -1$}
\put(150,10){$-1 \times i + (+1) \times j = -2$}
\end{picture}
\caption{Grading the points in the support along $(-1,+1)$.}
\label{fignormgrade}
\end{center}
\end{figure}

Given a vector $(u,v)$, $\langle \cdot , (u,v) \rangle$ ranks
the points $(i,j)$.
For $(u,v) = (1,1)$, we have the usual degree of $x^i y^j$.
So the direction of the tentacles are grading the points
in the support.
For example, in Figure~\ref{fignormgrade} we look
at the support in the direction $(-1,+1)$.

The degree of $x^i y^j$ in the direction $(u,v)$ is
the value of the inner product $\langle (i,j) , (u,v) \rangle$.  
In Maple we compute weighted degrees as follows:

\begin{center}
{\tt Groebner[WeightedDegree](f,[-1,+1],[x,y]);}
\end{center}

The connection between Gr\"obner bases and Newton polytopes
is the subject of~\cite{Stu96}.
This grading leads to homogeneous coordinates,
see~\cite{Cox03},

We arrive at a tropicalization of a polynomial via the normal fan
to the Newton polygon of the polynomial.

\begin{definition}  {\rm
Let $P$ be the Newton polygon of $f$.
The inner product is denoted by $\langle \cdot , \cdot \rangle$.
The {\em normal cone to a vertex $\bfp$ of $P$} is
\begin{equation}
  \{ \ \bfv \in \rr^2 \setminus \{ 0 \} \ | \ 
     \langle \bfp , \bfv \rangle = \min_{\bfq \in P} 
     \langle \bfq , \bfv \rangle \ \}.
\end{equation}
The {\em normal cone to an edge spanned by $\bfp_1$ and $\bfp_2$} is
\begin{equation}
  \{ \ \bfv \in \rr^2 \setminus \{ 0 \} \ | \ 
     \langle \bfp_1 , \bfv \rangle = \langle \bfp_2 , \bfv \rangle
     = \min_{\bfq \in P} \langle \bfq , \bfv \rangle \ \}.
\end{equation}
The {\em normal fan} of $P$ 
is the collection of all normal cones
to vertices and edges of $P$.
All normal cones to the edges of $P$ define
{\em a tropicalization of $f$}, denoted by ${\rm Trop}(f)$.
}
\end{definition}

We speak of {\em a} tropicalization ({\em a} instead of {\em the})
because different regular triangulations of the Newton polygon
will lead to different normal fans.  We can identify the tropicalization
by the secondary fan~\cite{GKZ94} but for our problem we do not need it.

\section{Tropisms}

The tropical view will lead to an efficient preprocessing stage
to determine whether two polynomials have a common factor.

\subsection{turning the varieties in a particular direction}

The answer to our original question
{\em Do two polynomials have a common factor?}
first depends on the relative position of the Newton polygons.
We compute tropicalizations of the polynomials and obtain
an efficient preprocessing step independent of the coefficients.

We first want to exclude the situations where there is no
common factor, already implied by the Newton polygons in
relative general position.  This is a direct consequence
of Bernshte\v{\i}n's second theorem~\cite{Ber75}.
For completeness, we state this theorem here for Newton polygons.

\begin{theorem}  \label{thebernshtein2nd}
Let $f$ and $g$ be two polynomials in $x$ and~$y$.
If ${\rm Trop}(f) \cap {\rm Trop}(g) = \emptyset$
then the system $f(x,y) = 0 = g(x,y)$ has no solutions at infinity. 
\end{theorem}

Now we can make Proposition~\ref{propamoebas} effective:

\begin{proposition}
If for two polynomials $f$ and $g$:
${\rm Trop}(f) \cap {\rm Trop}(g) = \emptyset$,
then $f$ and $g$ have no common factor.
\end{proposition}
\noindent {\em Proof.}  By Theorem~\ref{thebernshtein2nd},
${\rm Trop}(f) \cap {\rm Trop}(g) = \emptyset$
implies there is no common root at infinity.
But if $f$ and $g$ would have a common factor, they would
have a common root at infinity as well.~\qed

\begin{figure}[hbt]
\begin{center}
\begin{picture}(330,140)(0,0)
\put(0,10){\epsfig{figure=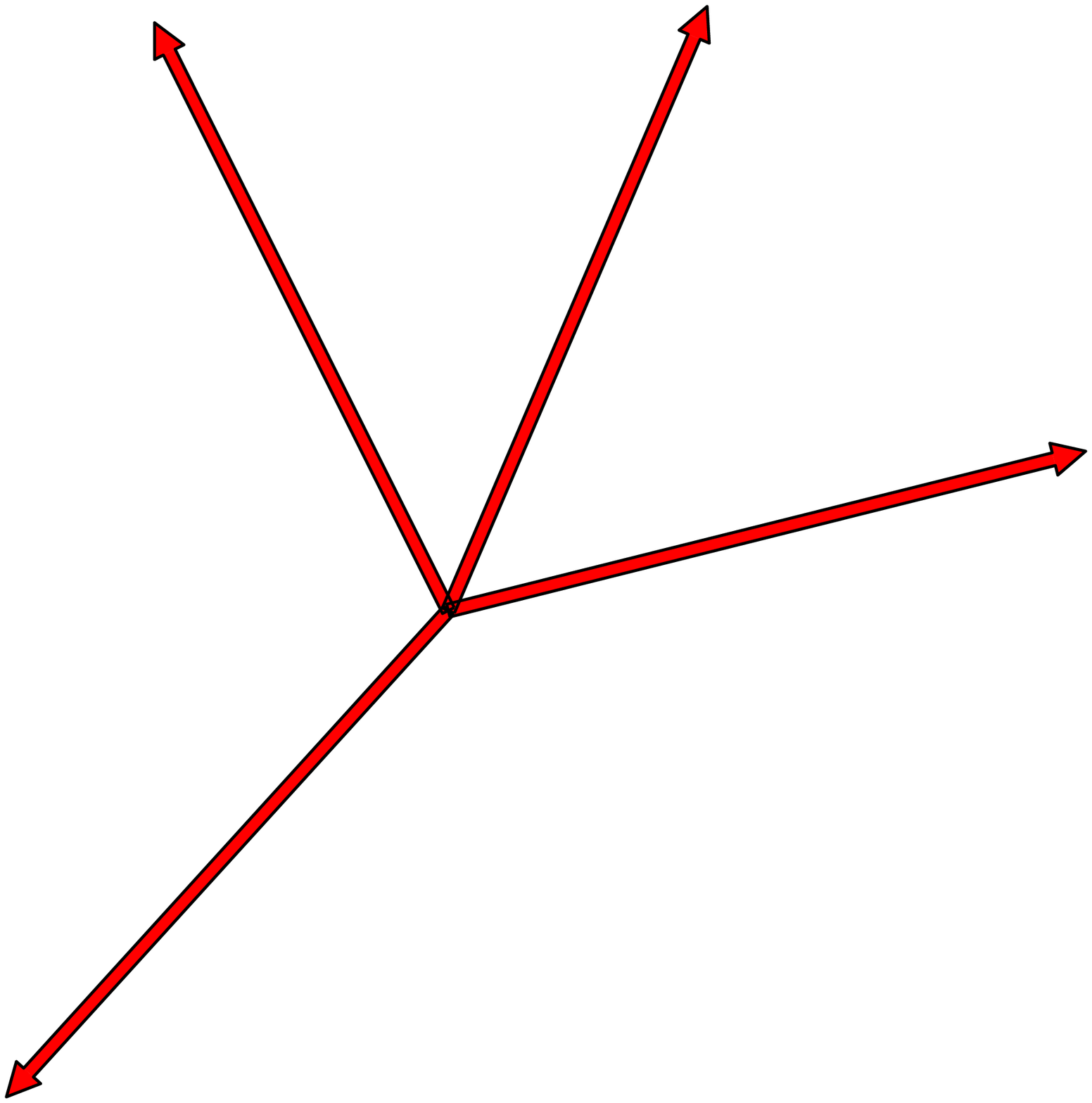,width = 4cm}}
\put(112,10){\epsfig{figure=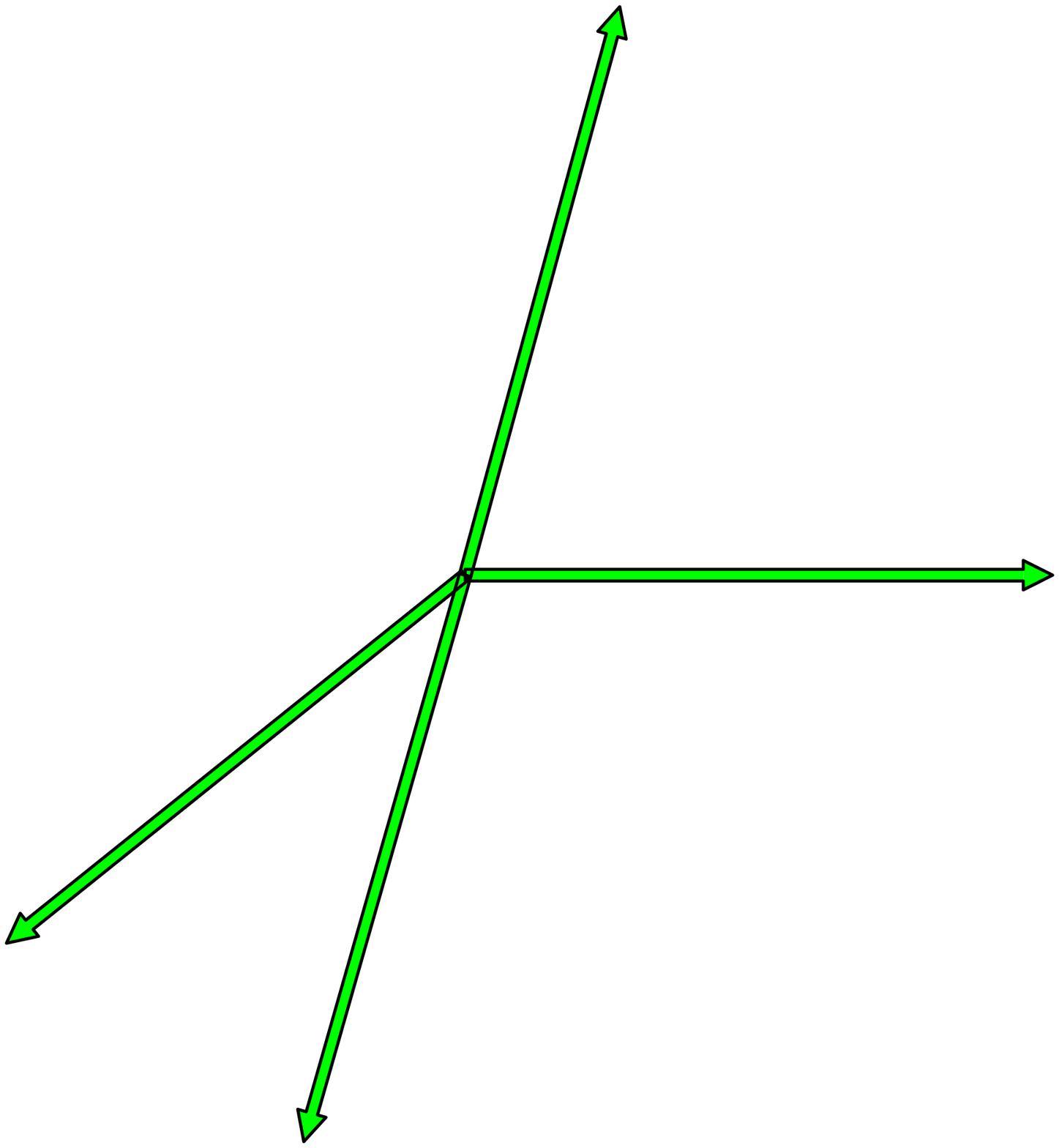,width = 4cm}}
\put(222,10){\epsfig{figure=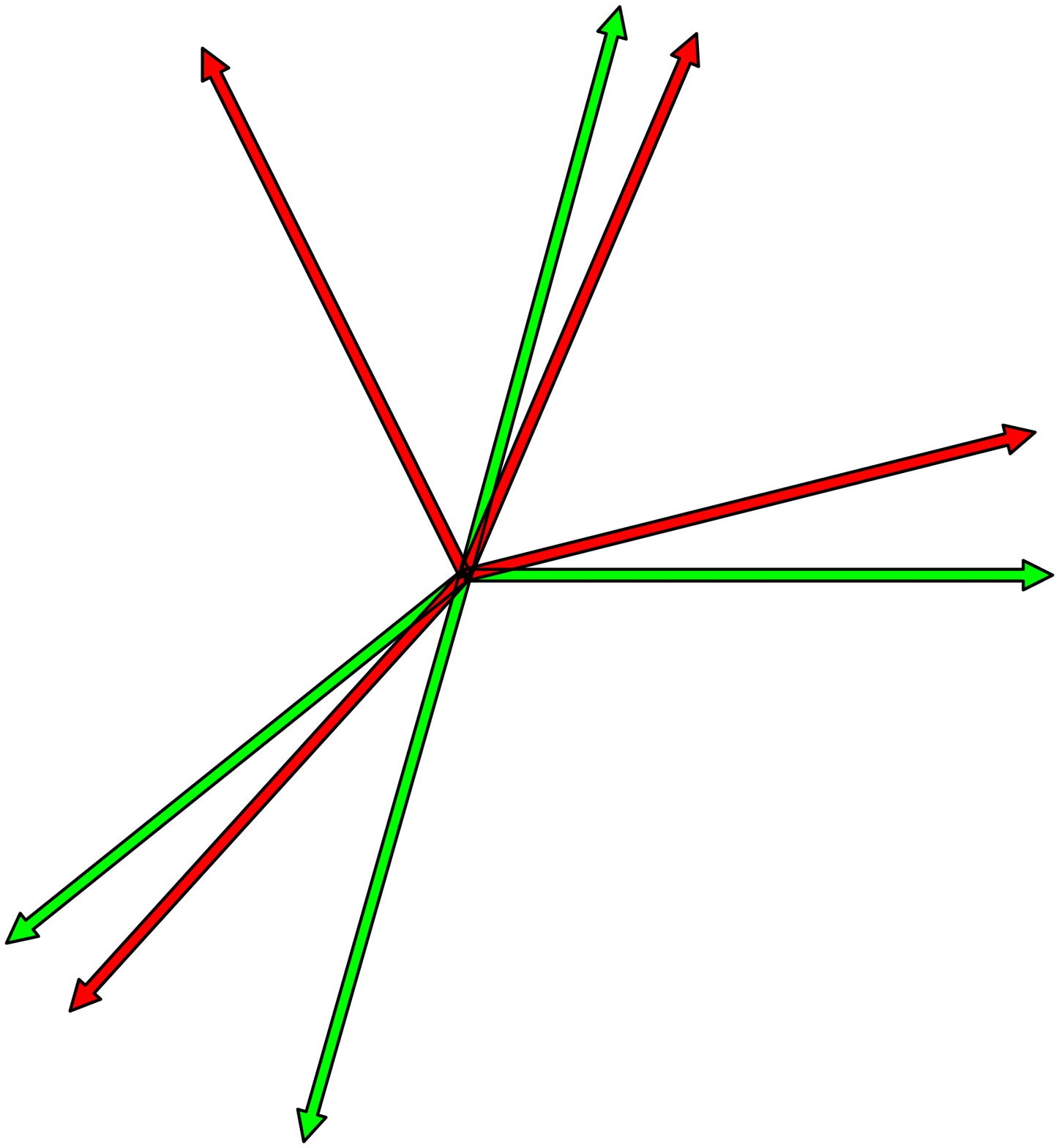,width = 4cm}}
\end{picture}
\caption{The first two pictures from the left represent the normal
fans of two polynomials.  By superposition of the fans at the far right
we see there are no common directions.  Therefore, for all nonzero
coefficients, the polynomials can have no common factor.}
\label{fignofacfan}
\end{center}
\end{figure}

Consider for example the tropicalizations of
two random polynomials of degree 15, shown in Figure~\ref{fignofacfan}.

\begin{figure}[hbt]
\begin{center}
\begin{picture}(330,140)(0,0)
\put(0,10){\epsfig{figure=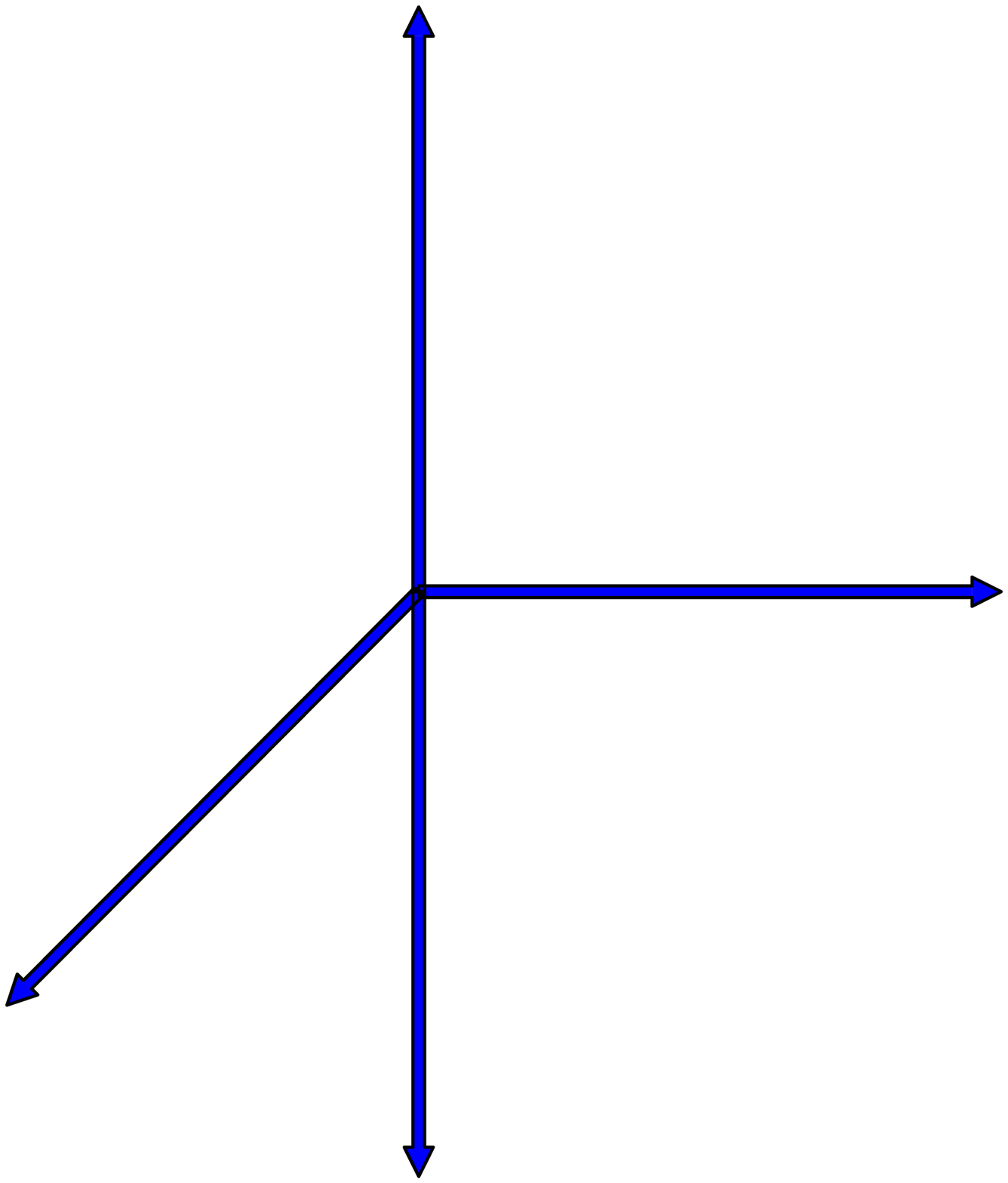,width = 4cm}}
\put(110,10){\epsfig{figure=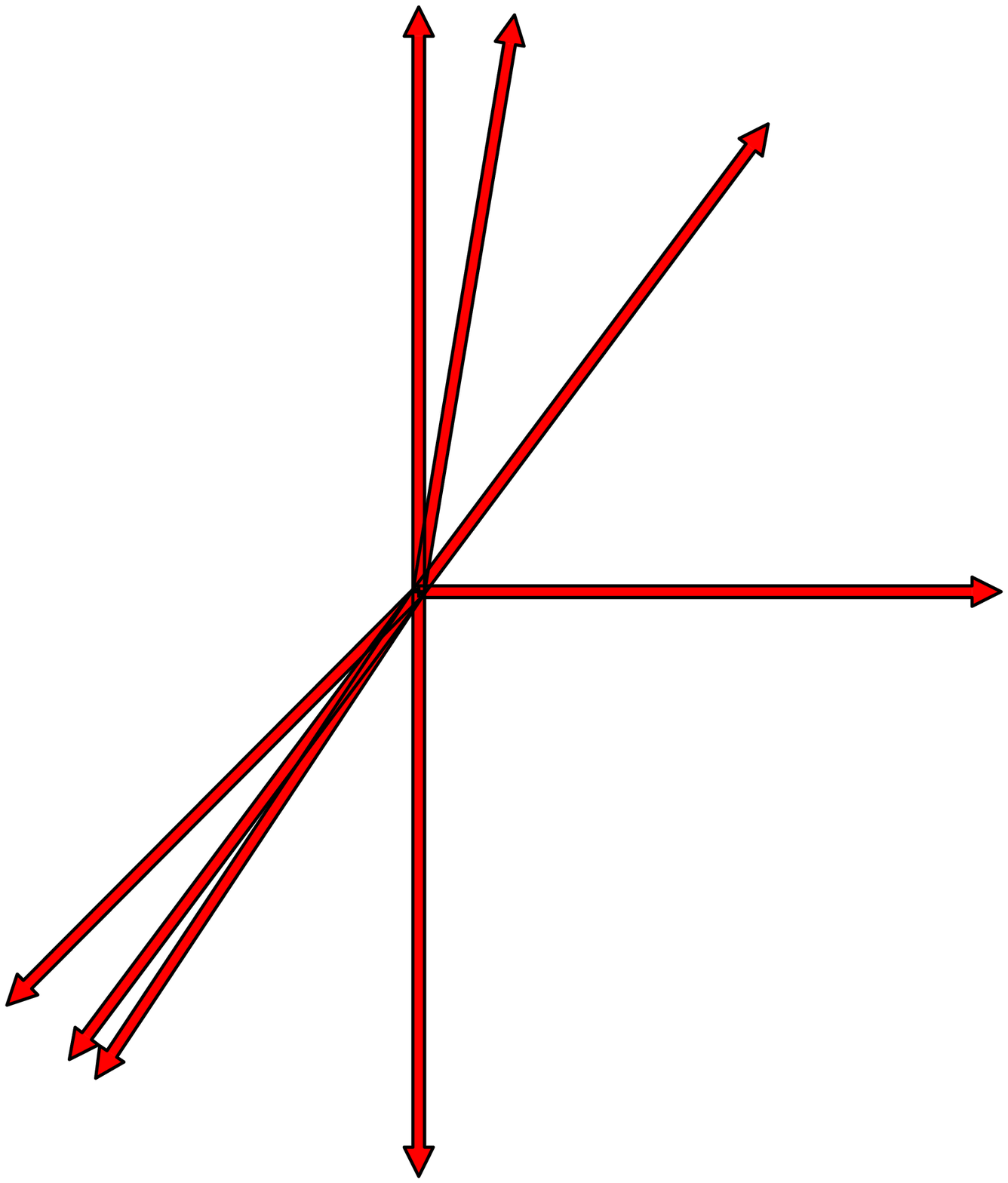,width = 4cm}}
\put(220,10){\epsfig{figure=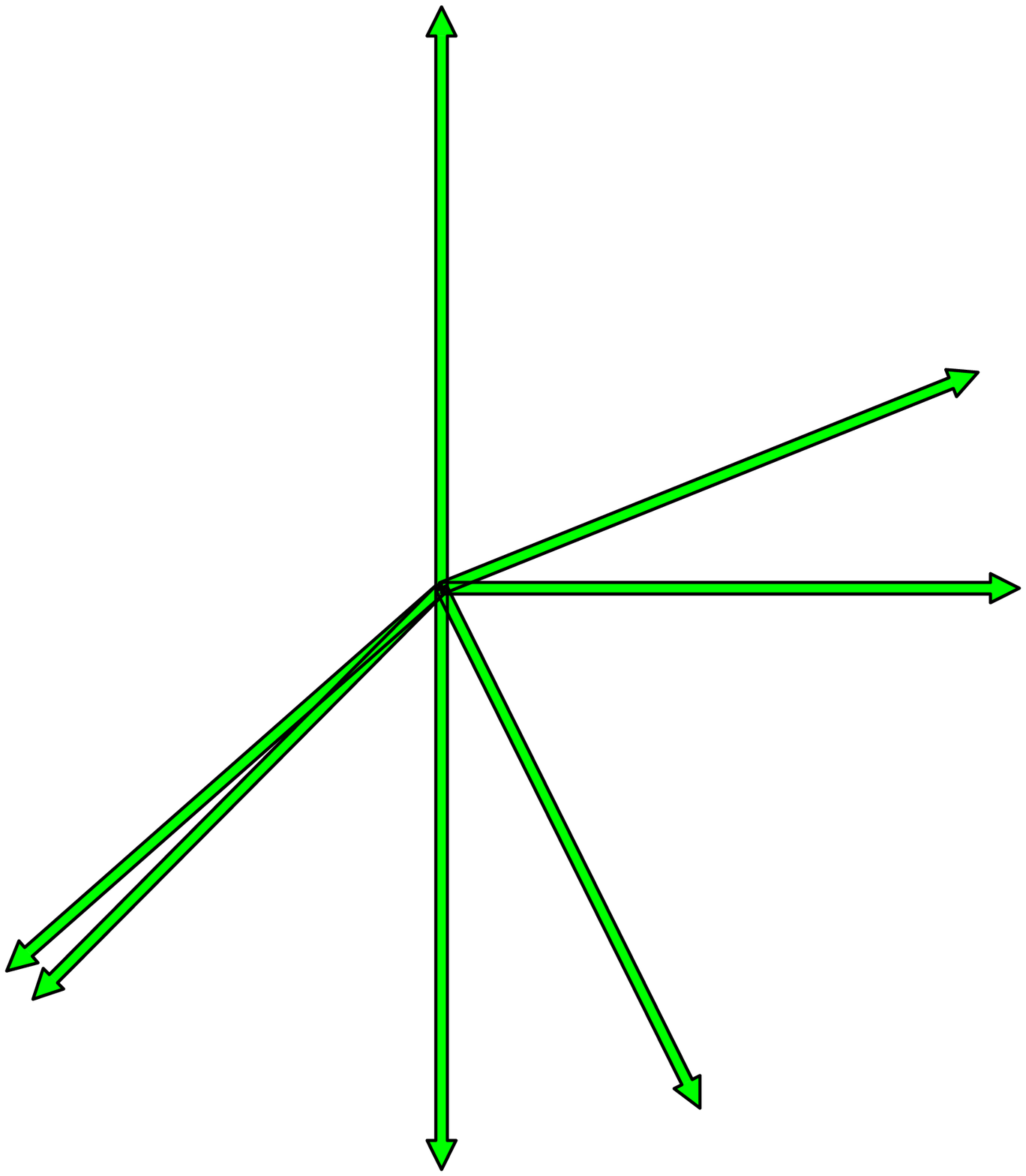,width = 4cm}}
\end{picture}
\caption{The normal fan at the left is the normal fan of the
factor common to two polynomials~$f$ and~$g$ whose normal fans
are displayed in the middle and at the right.  
We recognize the fan at the left as a part of the other fans.}
\label{figisfacfan}
\end{center}
\end{figure}

In our second example we generated a factor of degree 5 
and multiplied with two random polynomials of degree 10.
Tropicalization of the factor and the two polynomials are
shown in Figure~\ref{figisfacfan}.

A dictionary definition of a tropism is
{\em the turning of all or part of an organism in a particular 
direction in response to an external stimulus}.
Tropisms were introduced mathematically in 1980 
by Joseph Maurer~\cite{Mau80} who generalized Puiseux expansions 
for space curves.  We adapt his definition for use to our problem.

\begin{definition} {\rm
Let $P$ and $Q$ be Newton polygons of $f$ and $g$.
A {\em tropism} is a vector perpendicular
to one edge of $P$ and one edge of $Q$. }
\end{definition}

Using the general terminology of~\cite{Zie95},
tropisms correspond to the one dimensional cones in the common
refinement of the normal fans of the polygons.

\subsection{certificates for numerical computations}

Tropisms are important because they give a first exact
certificate for the existence of a common factor.
Selecting those monomials which span the edges picked out
by the tropism defines a polynomial system which admit a solution
in $(\cc^*)^2$.

\begin{definition} {\rm
Let $(u,v)$ be a direction vector.
Consider ${\displaystyle f = \sum_{(i,j) \in A} c_{i,j} x^i y^j}$.
The {\em initial form of $f$ in the direction $(u,v)$} is
\begin{equation}
{\rm in}_{(u,v)}(f)
= \sum_{\begin{array}{c}
           (i,j) \in A \\ \langle (i,j), (u,v) \rangle = m
        \end{array}} c_{i,j} x^i y^j,
\end{equation}
where $m = \min \{ \ \langle (i,j), (u,v) \rangle \ | \ (i,j) \in A \ \}$. }
\end{definition}
The direction $(u,v)$ is the normal vector to the line
$u i + v j = m$ which contains the edge of the Newton polygon
of $f$.  This edge is the Newton polygon of~${\rm in}_{(u,v)}(f)$.

The terminology of initial forms corresponds to the Gr\"obner
basics~\cite{Stu96}.  In~\cite{Roj03}, ${\rm in}_{(u,v)}(f)$ 
is called an initial term polynomial.
We call a tuple of initial forms {\em an initial form system}.
Initial form systems correspond to truncated systems
in~\cite{Bru00} and~\cite{Kaz99}.

\begin{figure}[hbt]
\begin{center}
\begin{picture}(330,130)(0,0)
\put(0,10){\epsfig{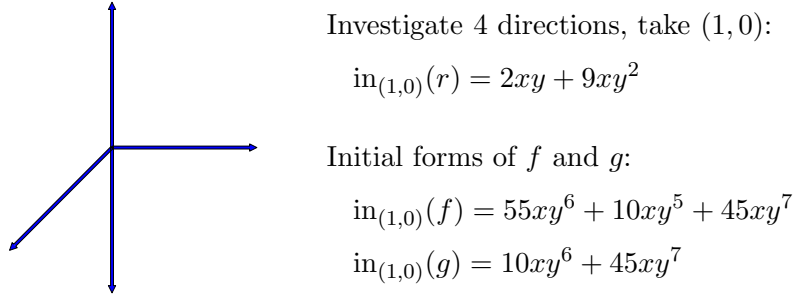}}
\put(130,110){Investigate 4 directions, take $(1,0)$:}
\put(140,90){${\rm in}_{(1,0)}(r) = 2 x y + 9 x y^2$}
\put(130,60){Initial forms of $f$ and $g$:}
\put(140,40){${\rm in}_{(1,0)}(f) = 55 x y^6 + 10 x y^5 + 45 x y^7$}
\put(140,20){${\rm in}_{(1,0)}(g) = 10 x y^6 + 45 x y^7$}
\end{picture}
\caption{The normal fan of the common factor and the 
initial form systems corresponding to the direction~(1,0).}
\label{figinforms}
\end{center}
\end{figure}

The factor common to the $f$ and $g$ generated above is
\begin{equation} \label{eqcomfac}
r := 2 x y + x^2 y + 9 x y^2 + 7 x^3 y + x^4 y + 9 x^3 y^2,
\end{equation}
In Figure~\ref{figinforms} we show the initial forms associated
with the two polynomials defined by the tropisms.

Because the tropism is a standard basis vector~(1,0),
the initial form system it determines consists of two
polynomials in one variable:
\begin{equation}
   \left\{
     \begin{array}{l}
       {\rm in}_{(1,0)}(f) = x \left( 5 y^5 (y+1)(2+9y) \right) = 0 \\
       {\rm in}_{(1,0)}(g) = x \left( 5 y^5 (2+9y) \right) = 0
     \end{array}
   \right.
\end{equation}
and then $y = -2/9$ represents the common root at infinity,
In general the common root at infinity will be 
an approximate root and with $\alpha$-theory~\cite{BCSS98}
we can bound the radius of convergence for Newton's method.
In addition to the first certificate, the exact tropism,
the common root at infinity is the second approximate certificate.

For general tropisms, not equal to basis vectors, 
we perform unimodular transformations
in the space of the exponents to reduce the initial form system
to a system of two polynomial in one variable.  
In~\cite{Bru00}, the coordinate transformations resulting from those
unimodular transformations are called power transformations and they
power up the field of ``Power Geometry''.

For example, investigating the direction $(-1,-1)$:
\begin{equation}
\left\{
  \begin{array}{rcl}
    {\rm in}_{(-1,-1)}(f) & = & 54 x^{13} y^2 + 6 x^{14} y
  \\ \vspace{-2mm} \\
    {\rm in}_{(-1,-1)}(g) & = & 72 x^9 y^{10} + 8 x^{10} y^9
  \end{array}
\right.
\end{equation}
We change coordinates using the unimodular matrix
$M = \left[ \begin{array}{rr} -1 & -1 \\ 0 & -1 \end{array} \right]$.

\begin{definition} \label{defcoordtrans}
{\rm 
For a tropism $(u,v)$ normalized so the
greatest common divisor ${\rm gcd}(u,v) = 1$,
the unimodular matrix~$M$
\begin{equation}
   M = \left[
     \begin{array}{cc}
         u & v \\ -l & k \\
     \end{array}
   \right],
   \quad {\rm gcd}(u,v) = 1 = k u + l v = \det(M)
\end{equation}
defines {\em the unimodular coordinate transformation}
$x = X^u Y^{-l}$ and $y = X^v Y^k$.
}
\end{definition}

Note that for a monomial $x^a y^b$, the coordinate transformation yields
\begin{equation}
(X^u Y^{-l})^a (X^v Y^k)^b = X^{a u + b v} Y^{-l a + k b}
 = X^{\langle (a,b) , (u,v) \rangle} Y^{-l a + k b},
\end{equation}
so after the coordinate transformation, the monomials in the
initial forms all have the same minimal degree in~$X$.

Continuing the example, we perform the change of coordinates:
\begin{equation}
\left\{
  \begin{array}{rcl}
    {\rm in}_{(-1,-1)}(f)(x = X^{-1}, y = X^{-1} Y^{-1})
 & = & ( 54 Y + 6 )/(X^{15} Y^{2})
  \\ \vspace{-2mm} \\
    {\rm in}_{(-1,-1)}(g)(x = X^{-1}, y = X^{-1} Y^{-1})
 & = & ( 72 Y + 8 )/(X^{19} Y^{10})
  \end{array}
\right.
\end{equation}
This change of coordinates reduces the initial form system
to a system of two polynomials in one variable.
For the example, $Y = -1/9$ represents the common root at infinity.
Going back to the original coordinates:
\begin{equation}
   \left\{
     \begin{array}{l}
        X =  t \\
        Y = -1/9 \\
     \end{array}
   \right.
   \quad 
   \left(
     \begin{array}{l}
         x = X^{-1} \\
         y = X^{-1} Y^{-1}
     \end{array}
   \right)
      \quad \Rightarrow \quad
   \left\{
      \begin{array}{l}
         x = t^{-1} \\
         y = -9 t^{-1}. \\
      \end{array}
   \right.
\end{equation}
As $t$ goes to~$0$ we have indeed a root going off to infinity.

For every tentacle of the common factor we can associate
a degree as follows.  Considering again
the common factor~$r$ from~(\ref{eqcomfac}),
the amoeba for $r$ has four tentacles, see Figure~\ref{figinforms},
reflected by its tropicalization
\begin{equation}
  {\rm Trop}(r) = \left\{ \ (1,0), (0,1), (-1,-1), (0,-1) \ \right\}.
\end{equation}
In Table~\ref{tabfacdegrees} we list the degrees associated
to each tentacle of the common factor.
We count the number of nonzero solutions of the initial forms,
after proper unimodular coordinate transformation.
To make the correspondence with the usual degree,
we consider only tropisms with nonnegative exponents
and ignore zero factors.

\begin{table}[h]
\begin{center}
\begin{tabular}{ccc}
  $(u,v)$ & ${\rm in}_{(u,v)}(r)$ & degree \\ \hline \\ \vspace{-7mm} \\
  $(1,0)$ & $2 x y + 9 x y^2$ &  1 \\ \vspace{-4mm} \\
  $(0,1)$ & $2 x y + x^2 y + 7 x^3 y + x^4 y$ & 3  \\ \vspace{-4mm} \\
  $(-1,-1)$ & $x^4 y + 9 x^3 y^2$ & 1 \\ \vspace{-4mm} \\
  $(0,-1)$ & $9 x y^2 + 9 x^3 y^2$ & 2 
\end{tabular}
\end{center}
\caption{Degrees associated to each vector in ${\rm Trop}(r)$.}
\label{tabfacdegrees}
\end{table}

\subsection{A preprocessing algorithm and its cost}

\begin{figure}[ht]
\begin{center}
\begin{picture}(200,200)(0,10)
\put(70,190){tropicalization}
\put(-10,156){1. compute tropisms}
\put(100,185){\vector(0,-1){15}} \put(100,160){\circle{8}}
\put(105,153){\vector(+1,-1){15}}
\put(125,135){no tropism}
\put(130,120){$\Rightarrow$ no root at $\infty$}

\put(-10,106){2. solve initial forms}
\put(100,150){\vector(0,-1){30}} \put(100,110){\circle{8}}
\put(105,103){\vector(+1,-1){15}}
\put(125,85){no root at $\infty$}
\put(130,70){$\Rightarrow$ no series}

\put(-10,56){3. compute 2nd term}
\put(100,100){\vector(0,-1){30}}  \put(100,60){\circle{8}}
\put(105,53){\vector(+1,-1){15}}
\put(125,35){no series}
\put(130,20){$\Rightarrow$ no factor}
\put(100,50){\vector(0,-1){30}}
\put(85,5){series}
\end{picture}
\end{center}
\label{figstages}
\caption{A staggered approach to computing a common factor.}
\end{figure}
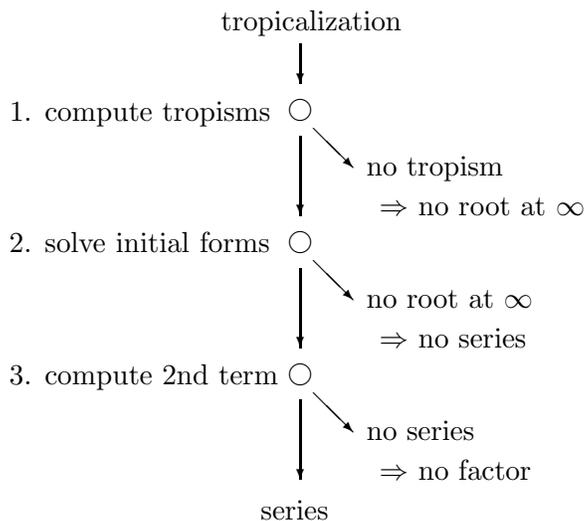

That two polynomials with approximate coefficients have
a common factor does not happen that often.  Therefore,
it is important to be able to decide quickly in case there
is no common factor.  The stages in a preprocessing algorithm
are sketched in Figure~\ref{figstages}.

In Figure~\ref{figstages} we distinguish three computational steps.
We will address the cost of the first two steps in the following
propositions.

\begin{proposition} \label{propcost1}
Let $f$ and $g$ be two polynomials given by respectively
$n$ and $m$ monomials.
The cost of computing tropisms ${\rm Trop}(f) \cap {\rm Trop}(g)$
is $O(n\log(n)) + O(m\log(m))$.
\end{proposition}

\noindent {\em Proof.}
It takes $O(n \log(n))$ operations for computing
a tropicalization ${\rm Trop}(f)$ because computing the convex
hull of a set of $n$ points amounts to sorting the points
in the support.  Likewise, computing ${\rm Trop}(g)$ takes
$O(m \log(n))$ operations.
Merging sorted lists of normals to find the tropisms
in ${\rm Trop}(f) \cap {\rm Trop}(g)$ takes linear time
in the length of the lists.~\qed

This preprocessing step has the lowest complexity
and as the algorithm operates only on the exponents
the outcome is exact.  The absence of tropisms is
an exact certificate that there is no common factor,
for any nonzero choice of the coefficients of the polynomials.

In case we have tropisms, we have to compute roots at infinity.
The cost of the second preprocessing stage is as follows.

\begin{proposition}  \label{propcost2}
Let $f$ and $g$ be two polynomials
given by respectively $n$ and $m$ monomials.
For every tropism $t \in {\rm Trop}(f) \cap {\rm Trop}(g)$
it takes at most $O((n + m)^3)$ operations to find a common
solution in $(\cc^*)^2$ to the initial form system defined by~$\bfv$.
\end{proposition}

\noindent {\em Proof.}  For a tropism~$\bfv$, 
we solve the initial form system.  In particular,
an initial root $\bfz$ satisfies

\begin{equation}
\left\{
   \begin{array}{c}
      {\rm in}_\bfv(f)(\bfz) = 0 \\
      {\rm in}_\bfv(g)(\bfz) = 0 \\
   \end{array}
\right. \quad \bfz \in (\cc^*)^2.
\end{equation}
We perform a unimodular transformation so the tropism
we consider is a unit vector, (1,0) or (0,1).
This implies that the two equations in the initial form system
are defined by two polynomials in one variable.
To decide whether two polynomials in one variable admit a common
solution we determine the rank of the Sylvester matrix.
Using singular value decomposition, the cost of this rank
determination is cubic in the size of the matrix.~\qed

Even as the cost estimates in the propositions are conservative,
they give a good polynomial complexity.  Actually, in the best case,
the initial forms are supported on two points only and instead of
a rank determination, we can just take primitive roots.
The cost estimates of Proposition~\ref{propcost2} cover the very
worst situation where the Newton polygons are triangles and one
of the edges contains all exponent vectors except for one.

For numerical calculations, it is important to note that in this
preprocessing stage, only the coefficients at the edges are involved.
If the coefficients are badly scaled, then coefficients with
monomials in the interior of the Newton polygons will not cause
difficulties in this preprocessing stage.

For the complexity in the proof of the second proposition we used
the ubiquitous singular value decomposition
but for practical purposes rank revealing algorithms~\cite{LZ05}
have a lower cost.
The accurate location of the root of the initial form systems may
look complication in case this root is multiple.  However, because
the initial form systems consists of univariate equations,
the methods of~\cite{Zen05} will give satisfactory answers.

\section{Germs}

Once we have a tropisms and an initial root at infinity,
we start growing the Puiseux series for the common factor.

\subsection{how the amoeba grows from infinity}

We use the roots at infinity to grow the tentacles
of the common factor.  But first we must decide
whether the roots at infinity are isolated or not.

The canonical form of a fractional power series 
Puiseux series of a curve is defined next.

\begin{definition} {\rm
Consider the curve defined by $r(x,y) = 0$.
A {\em fractional power series} has the form
\begin{equation}
  \left\{
     \begin{array}{ccl}
        x & = & z_1 t^u \\
        y & = & z_2 t^v (1 + O(t))
     \end{array}
  \right.
  \quad (z_1,z_2) \in (\cc^*)^2.
\end{equation} }
\end{definition}
The leading exponents $(u,v)$ are the tropisms.
The leading coefficients $(z_1,z_2)$ satisfy
the initial form system
\begin{equation}
   \left\{
      \begin{array}{c}
         {\rm in}_{(u,v)}(f)(z_1,z_2) = 0 \\
         {\rm in}_{(u,v)}(g)(z_1,z_2) = 0 \\
      \end{array}
   \right.
\end{equation}
where the initial forms are taken from the equations $f$
and $g$ which define the common factor~$r$.

Continuing our example, we extend the solution at infinity,
defined by the initial form system for the first tropism~(1,0).
Because the tropism is a standard basis vector, the Maple
command {\tt sort(\{ f , g \}, plex, ascending)}
will show that the leading terms of the polynomials $f$ and $g$
are indeed ${\rm in}_{(1,0)}(f)$ and ${\rm in}_{(1,0)}(g)$:
\begin{equation}
\left\{
  \begin{array}{l}
     f = 10 x y^5 + 45 x y^7 + 55 x y^6 + x^2 ( \mbox{ 30 other terms } ) 
\\ \vspace{-2mm} \\
     g = 45 x y^7 + 10 x y^6 + x^2 ( \mbox{ 34 other terms } )
  \end{array}
\right.
\end{equation}

Let $f_1 = f/x$ and $g_1 = g/x$, 
then $z = -2/9$ is solution at infinity.

\begin{equation}
  \left\{
     \begin{array}{l}
        x = t^{1} \\
        y = -\frac{2}{9} t^{0} + C t (1 + O(t)), \quad c \in \cc^*.
     \end{array}
  \right.
\end{equation}
A nonzero value for $C$ will give the third certificate
of a common factor.
Useful Maple commands to compute the power series are

\begin{verbatim}
zt := x = t, y = -2/9 + C*t;
f1z := subs(zt,f1): g1z := subs(zt,g1):
c1 := coeff(f1z,t,1); c2 := coeff(g1z,t,1);
\end{verbatim}
The constraints on the coefficient~$C$ we obtain are
\begin{equation}
  \left\{
     \begin{array}{l}
        c1 = -\frac{1120}{531441} - \frac{1120}{59049} C = 0 \\ 
  \vspace{-2mm}  \\
        c2 = -\frac{320}{59049} - \frac{320}{531441} C = 0
     \end{array}
  \right.
\end{equation}
Notice that the second coefficient~$C$ of the Puiseux series expansion
again must satisfy an overdetermined system.
Solving both equations for $C$ gives $C = -1/9$.
\begin{equation}
  \left\{
    \begin{array}{l}
       x = t \\
       y = - \frac{2}{9} - \frac{1}{9} t (1 + O(t)).
    \end{array}
  \right.
\end{equation}
Substituting $x = t, y = -2/9 - t/9$ into $f_1$ and $g_1$
gives $O(t^2)$.  The second term of the Puiseux series is
the third and last certificate for a common factor.

In general, the next term in the Puiseux series expansion might
have a degree higher than one, or they might not exist a second
term at all in case the solution at infinity is isolated.
There is an explicit condition on the exponent of the second
term in the Puiseux series expansion as in Proposition~\ref{propsecexp}.

\begin{proposition} \label{propsecexp}
Let $f$ and $g$ be two polynomials in two variables
with $(u,v)$ a tropism for which there is a root of the
corresponding initial form system.
After a unimodular coordination defined by the tropism $(u,v)$
the Puiseux series takes the standard form
\begin{equation} \label{eqstandard}
  \left\{
     \begin{array}{l}
        x = t^d \\
        y = c_0 + c_1 t^w
     \end{array}
  \right.
  \quad {\rm with} ~ d = \gcd(u,v)
\end{equation}
and $c_0$ is a nonzero solution of the initial form system:
\begin{equation}
   \left\{
      \begin{array}{c}
         p(c_0) = 0 \\ q(c_0) = 0
      \end{array}
   \right.
   \quad {\rm where} \quad
   \left\{
      \begin{array}{l}
         {\rm in}_{(1,0)} f(x,y) =  p(y) \\
         {\rm in}_{(1,0)} g(x,y) =  q(y).
      \end{array}
   \right.
\end{equation}
Assuming that
\begin{equation} \label{eqassuming}
   \left\{
      \begin{array}{lcr}
         f(x,y) = p(y) + x^k P(x,y), & k > 0,
       & P(x,y) = p_{00} + p_{10} x + p_{01} y + \cdots\\
         g(x,y) = q(y) + x^l Q(x,y), & l > 0,
       & Q(x,y) = q_{00} + q_{10} x + q_{01} y + \cdots\\
      \end{array}
   \right.
\end{equation}
so that we may write
$p$ and $q$ as follows:
\begin{equation}
    \begin{array}{ccccccccc}
        p \!\! & = & \!\!\! \gamma_0 \!\!\! 
          & + & \!\!\! \gamma_1 y^{a_1} \!\!\! 
          & + & \!\!\! \gamma_2 y^{a_2} \!\!\! 
          & + & \!\!\! \cdots \\
          &   & 0 & < & a_1 & < & a_2 & < & \!\!\! \cdots
    \end{array}
\quad {\rm and} \quad
    \begin{array}{ccccccccc}
        q \!\! & = & \!\!\! \delta_0 \!\!\! 
          & + & \!\!\! \delta_1 y^{b_1} \!\!\! 
          & + & \!\!\! \delta_2 y^{b_2} \!\!\! 
          & + & \!\!\! \cdots \\
          &   & 0 & < & b_1 & < & b_2 & < & \!\!\! \cdots
    \end{array}
\end{equation}
then the condition on the exponent $w$ such that there may
exist a nonzero $c_1$ is
\begin{equation} \label{eqcondexp}
   w = \frac{k d}{a_1} = \frac{l d}{b_1}.
\end{equation}
If $\frac{k \times d}{a_1} \not= \frac{l \times d}{b_1}$,
then there is no second term in the Puiseux series expansion
and the root of the initial form system is isolated.
\end{proposition}

\noindent {\em Proof.}  We substitute the series~(\ref{eqstandard})
into the assumed form of the system~(\ref{eqassuming}) and look for
the lowest $w$ so that the coefficient with the monomial of the lowest
power of~$t$ vanishes.  For the polynomials in the initial form system
the substitution yields
\begin{equation}
   p(c_0 + c_1 t^w) = \alpha_1 c_1 t^{a_1 w} + h.o.t. \quad {\rm and} \quad
   q(c_0 + c_1 t^w) = \beta_1 c_1 t^{b_1 w} + h.o.t.
\end{equation}
and for the other terms $x^k P(x,y)$ and $x^l Q(x,y)$ we find
\begin{equation}
   t^{kd} ( \alpha_2 + O(t)), \alpha_2 = p_{00} + p_{01} c_0
    \quad {\rm and} \quad
   t^{ld} ( \beta_2 + O(t)), \beta_2 = q_{00} + q_{01} c_0.
\end{equation}
Collecting the lowest degree terms in~$t$ we find the system
\begin{equation} \label{eqlowestdeg}
  \left\{
     \begin{array}{c}
        \alpha_1 c_1 t^{a_1 w} + \alpha_2 t^{kd} = 0 \\
        \beta_1 c_1 t^{b_1 w} + \beta_2 t^{ld} = 0. \\
     \end{array}
  \right.
\end{equation}
For $c_1$ to exists as a nonzero solution of a linear system
annihilating the lowest degree terms in~$t$, the exponents of~$t$ 
in the system~(\ref{eqlowestdeg}) must be equal in both equations:
\begin{equation}
   a_1 w = kd \quad {\rm and} \quad b_1 w = l d.
\end{equation}
Eliminating $w$ leads to the condition~(\ref{eqcondexp}).~\qed

The assumption in~(\ref{eqassuming}) of Proposition~\ref{propsecexp}
makes abstraction of the monomial shifts that are usually required
when looking for nonzero solutions.  So the assumption~(\ref{eqassuming})
will not hold in general and monomial shifts are needed.
But this will still lead to an explicit condition on the exponents
purely derived from the exponents of the polynomials $f$ and $g$
for a second term in the Puiseux series to exist, leading to a
linear system of two equations in the one coefficient $c_1$.

In case the condition of Proposition~\ref{propsecexp} is satisfied
and the linear system admits a nonzero solution for~$c_1$, then
the exponent~$w$ and coefficient~$c_1$ constitute respectively
an exact and an approximate certificate for the existence of a
common factor for the two polynomials~$f$ and~$g$.

\subsection{regions of convergence of Puiseux series}

Neighborhoods where Puiseux series converge are called germs.
Following~\cite{dJP00}, we get the definition for a germ:

\begin{definition} {\rm Given a power series expansion at a point,
a {\em germ} is an equivalence class 
of open neighborhoods of the point where the series converge. }
\end{definition}

Via substitution as in the proof of Proposition~\ref{propsecexp}
we may compute more terms in the Puiseux series.  In each step,
the exponent in~$t$ will increase by at least one.
In each step, also more and more approximate coefficients of the 
given polynomials will be used.  The accuracy of the given input 
coefficients will impose a natural bound on the order for which
it still makes sense to extend the series.

The unimodular coordinate transformations play a very important
role also in the accurate evaluation of polynomials~\cite{DDH06}.
As the size of arguments of the polynomial functions grows,
and as the direction of the growth points along the direction
of a tentacle of the amoeba, monomials on the faces perpendicular
to that direction become dominant.  A weighted projective
transformation as in~\cite{Ver00} will rescale the problem
of evaluating a high degree polynomial with approximate coefficients
near a root.

\section{Implementation Aspects}

For efficient implementation of the algorithm, the data structures
used to represent the polynomials consist of a list of exponent vectors
and a coefficient table.  More precisely, to represent a polynomial~$f$
denoted as
\begin{equation}
  f(\bfx) = \sum_{\bfa \in A} c_\bfa \bfx^\bfa
  \quad c_\bfa \in \cc^*, \bfx^\bfa = x_1^{a_1} x_2^{a_2}
\end{equation}
we use a list to represent the support~$A$ and a lookup table $C[A]$
for the coefficients.
The indices of the lookup table~$C_A$ are the exponent vectors~$\bfa \in A$.
In Maple's index notation: $C[\bfa] = c_\bfa$.

Separating the support from the coefficient allows an
efficient execution of change of monomial orders.
If $n = \#A$, then monomial orders on~$f$ are stored via permutations of
the first $n$ natural numbers.  The separation also gives an efficient
way to change coordinates, i.e.: we apply the unimodular coordinate 
transformation only on~$A$.  For a unimodular matrix~$M$:
\begin{equation}
   MA = \{ \ M \bfa \ | \ \bfa \in A \ \}.
\end{equation}
Abusing notation, for $z \in \cc^*$: $Mz$ denotes the value for $Y$
after applying the coordinate transformation as in
Definition~\ref{defcoordtrans}.

The input polynomials $f$ and $g$ with respective supports $A_f$ and $A_g$
are then represented by two tuples: $(A_f,C[A_f])$ and $(A_g,C[A_g])$.
The preprocessing algorithms consists of two stage.
In the first stage, Algorithm~\ref{algfirststage} computes the tropisms 
and the roots of the corresponding initial form systems.  
If the sets of roots are not empty, the exponent and coefficients of 
the second term in the Puiseux expansions are computed 
by Algorithm~\ref{algsecondstage} in the second stage.
We define the specifications of the algorithms below.

\begin{algorithm} \label{algfirststage}
{\rm Tropisms and Initial Roots

\begin{tabular}{rcl}
  Input & : & $(A_f,C[A_f])$ and $(A_g,C[A_g])$. \\
 Output & : & $T = \{ \ (u,v) \in \zz^2 \setminus (0,0) \ | \
                        (u,v) ~\mbox{ is tropism } \}$, \\
        &   & $R[T] = \{ \ \{ \ z \in \cc^* \ | \
              {\rm in}_{(u,v)}(f)(Mz) = 0, 
              {\rm in}_{(u,v)}(g)(Mz) = 0 \ \} \ | \ (u,v) \in T \ \}$.
\end{tabular}
}
\end{algorithm}
Every tropism in~$T$ defines a set of roots (possibly empty)
of the corresponding initial form system, after application
of the unimodular coordinate transformation~$M$.
The cost of Algorithm~\ref{algfirststage} is estimated by
Proposition~\ref{propcost1} and Proposition~\ref{propcost2}.

\begin{algorithm} \label{algsecondstage}
{\rm Second Term of Puiseux Expansion

\begin{tabular}{rcl}
  Input & : & $(A_f,C[A_f])$, $(A_g,C[A_g])$, $T$, and $R[T]$. \\
 Output & : & $W[R[T]] = \{ \ (c,w) \in \cc^* \times \nn^+ \ | \ 
                z \in Z \in R[T] \ \}$.
\end{tabular}
}
\end{algorithm}
The elements of the set $W[R[T]]$ define the second term of the
Puiseux series expansion.  In particular, 
for every $(c,w) \in W[R[T]]$:
\begin{equation}
  \left\{
    \begin{array}{l}
       X = t^d \\
       Y = z + c t^w
    \end{array}
  \right.
\end{equation}
where $(X,Y)$ are the new coordinates after applying
the transformation of Definition~\ref{defcoordtrans}.
Conditions on the existence of the exponent~$w$
are given in Proposition~\ref{propsecexp}.

\section{Conclusions and Extensions}

Like Maple, tropical algebraic geometry is language.
Sentences like
{\em tropisms give the germs to grow
     the tentacles of the common amoeba}
express efficient preprocessing stages to detect and compute
common factors of two polynomials with approximate coefficients.
Seeing the problem as a system of two polynomial equations
in two variables, the algorithm is a polyhedral method to
find algebraic curves.

Among the extensions we consider for future developments
are algorithms to handle singularities numerically and
polyhedral methods for space curves.

\bibliographystyle{plain}

\end{document}